\documentclass[letterpaper,twoside,10pt]{ieeeconf}
\usepackage{generic}
\usepackage{cite}
\usepackage{amsmath,amssymb,amsfonts}
\usepackage{algorithmic}
\usepackage{graphicx}
\usepackage{tikz}
\usepackage{pgfplots}
\pgfplotsset{compat=1.18}
\usepackage{empheq}
\usepackage{comment}
\newtheorem{theorem}{Theorem}
\newtheorem{proposition}{Proposition}
\newtheorem{corollary}{Corollary}
\newtheorem{property}{Property}
\newtheorem{remark}{Remark}
\newtheorem{example}{Example}
\newtheorem{definition}{Definition}
\newtheorem{ass}{Assumption}
\newtheorem{lem}{Lemma}

\newenvironment{assumption}{\begin{ass}}{\hfill $\bullet$ \end{ass}}
\newenvironment{definitionn}{\begin{definition}}{\hfill $\bullet$ \end{definition}}

\newenvironment{Remark}{\begin{remark}}{\hfill $\bullet$ \end{remark}}
\newenvironment{Theorem}{\begin{theorem}}{\hfill $\square$ \end{theorem}}
\newenvironment{Proposition}{\begin{proposition}}{\hfill $\square$ \end{proposition}}

\newenvironment{Property}{\begin{property}}{\hfill $\bullet$ \end{property}}

\usepackage{algorithm,algorithmic}
\usepackage{hyperref}

\hypersetup{hidelinks=true}
\usepackage{textcomp}
\def\BibTeX{{\rm B\kern-.05em{\sc i\kern-.025em b}\kern-.08em
    T\kern-.1667em\lower.7ex\hbox{E}\kern-.125emX}}
\begin{document}
\title{Boundary Stabilization of Quasilinear Parabolic PDEs \\
that Blow Up in Open Loop for Arbitrarily Small Initial Conditions}
\author{M. C. Belhadjoudja, M. Maghenem, E. Witrant, M. Krstic
\thanks{M. C. Belhadjoudja and M. Maghenem are with Universit\'e Grenoble Alpes, CNRS, Grenoble-INP, GIPSA-lab, F-38000, Grenoble, France (e-mail: mohamed.belhadjoudja,mohamed.maghenem@gipsa-lab.fr).}
\thanks{E. Witrant is with Universit\'e Grenoble Alpes, CNRS, Grenoble-INP, GIPSA-lab, F-38000, Grenoble, France, and the Departement of Mechanical Engineering, Dalhousie University, Halifax B3H 4R2, Nova Scotia, Canada (e-mail: emmanuel.witrant@gipsa-lab.fr).}
\thanks{M. Krstic is with the Department of Mechanical and Aerospace Engineering, University of California San Diego, 92093 San Diego, USA (e-mail:
krstic@ucsd.edu).}}

\maketitle

\pagenumbering{gobble}

\begin{abstract}
We propose a novel framework for stabilization, with an estimate of the region of attraction, of quasilinear parabolic partial differential equations (PDEs) that exhibit finite-time blow-up phenomena when null boundary inputs are imposed. Using Neumann-type boundary controllers, which are cubic polynomials in boundary measurements, we ensure $L^2$ exponential stability of the origin with an estimate of the region of attraction, boundedness and exponential decay towards zero of the state's max norm, well-posedness, as well as positivity of solutions starting from positive initial conditions. Unlike existing methods, our approach handles nonlinear state-dependent diffusion, convection, and reaction terms. In many cases, our estimate of the size of the region of attraction is shown to expand unboundedly as diffusion increases. Our controllers can be implemented as Neumann, Dirichlet, or mixed-type boundary conditions. Numerical simulations validate the effectiveness of our approach in preventing finite-time blow up.
\end{abstract}

\section{introduction}
Quasilinear parabolic PDEs exhibiting finite-time blow-up phenomena appear across diverse engineering domains, including combustion theory \cite{combustion1,combustion2}, chemical kinetics \cite{chemical1,chemical2}, and plasma physics \cite{plasma}. These mathematical singularities model catastrophic events - such as thermal runaway in lithium-ion batteries leading to explosion, reactor core meltdowns, or uncontrolled chain reactions - where the system's state reaches extreme values in finite time. This has naturally motivated the control community to study these equations, aiming to prevent finite-time blow up through suitable control strategies. Unfortunately, \cite[Theorem 1.1]{lack2} established a key limitation: even for simple instances of the semilinear heat equation, global stabilization cannot be achieved with neither boundary nor in-domain control. Worse, for sufficiently large initial conditions, finite-time blow up is inevitable, regardless of the choice of control inputs. Thus, the stabilization of quasilinear parabolic PDEs comes in the form of an estimate of the region of attraction. Nevertheless, estimating the region of attraction remains crucial, as finite-time blow up can occur in open loop even for arbitrarily small initial conditions.

While boundary control of parabolic PDEs has been thoroughly investigated in the literature \cite{burgers,backstepping,lyapunov2}, to our knowledge, the only studies addressing the problem of stabilization of parabolic equations that may exhibit finite-time blow up phenomena are \cite{volterra1} and \cite{volterra2}, but are limited to the specific subclass of semilinear diffusion-reaction equations with Volterra-type nonlinearities, and use nonlinear backstepping. This approach cannot, however, handle quasilinear parabolic PDEs with nonlinear state-dependent diffusion, or equations with nonlinear convection. Furthermore, the positivity of solutions starting from positive initial conditions is not preserved. That being said, \cite{shock} addresses the regulation, via backstepping, of shock-like equilibria for Burgers' equation, which includes quadratic convection, and for which the open-loop error PDE exhibits finite-time blow up. However, the approach in \cite{shock} relies on the Cole-Hopf transform to linearize convection-a change of variables specific to Burgers’ equation. We mention also \cite{lyapunov1}, which achieves global stabilization for a class of semilinear diffusion-reaction equations with polynomial reaction terms, but under the assumption that the nonlinearity is stabilizing when the state is large.

In this paper, we develop a novel framework for stabilization of quasilinear parabolic PDEs with an estimate of the region of attraction. We design Neumann-type boundary controllers as cubic polynomials of the state values at the boundaries, that guarantee: (i) exponential stability of the origin in the $L^2$ norm with an estimate of the region of attraction, (ii) boundedness and exponential decay towards zero of the max norm of the state, (iii) existence and uniqueness of complete classical solutions, and (iv) positivity of solutions, i.e. the solutions are nonnegative if the initial conditions are nonnegative. Additionally, our Neumann-type boundary controllers can be rewritten as Dirichlet-type or mixed Dirichlet-Neumann-type controllers by solving some cubic equations using Cardano's root formula. In many cases, the region of attraction--as well as the decay rate--are shown to grow unboundedly as the diffusion increases. In contrast, when the boundary inputs are set to zero, finite-time blow up can occur for an arbitrarily small initial condition and an arbitrarily large diffusion coefficient. Our controllers are inspired by \cite{burgers}, the first work to ever use cubic Neumann feedback, but in a fundamentally different context from the one considered here. Indeed, \cite{burgers} considers the viscous Burgers' equation, which does not exhibit finite-time blow-up phenomena.

The paper is organized as follows. In Section \ref{system description}, we formulate our problem. In Section \ref{main results}, we present our main results. Finally, in Section \ref{simu}, we illustrate our results via numerical simulations, before we finalize the paper with a conclusion and research perspectives.

\textit{Notation.} Given $v : [0,1] \to \mathbb{R}$, we let $|v|_{\infty} := \max_{x \in [0,1]}|v(x)|$. Similarly, given $T \in (0,+\infty]$,  
$u : [0,1] \times [0,T) \to \mathbb{R}$, and $t \in [0,T)$, 
 we let $|u(\cdot,t)|_{\infty}:=\max_{x \in [0,1]}|u(x,t)|$. We also write $|u|_{\infty}$ to mean the map $t \mapsto |u(\cdot,t)|_{\infty}$.
Furthermore, we denote by $u_x$ the partial derivative of $u$ with respect to $x$, $u_{xx}$ the second partial derivative of $u$ with respect to $x$, and $u_t$ the partial derivative of $u$ with respect to $t$. Given $\beta \in (0,1)$, we denote by $\mathcal{H}^{2+\beta}[0,1]$ 
the space of twice continuously-differentiable maps $v: [0,1] \rightarrow \mathbb{R}$ such that $v''$, the second-order derivative of $v$, is Hölder continuous with exponent $\beta$, i.e., there exists $C>0$ such that 
$ |v''(x_1)-v''(x_2)| \leq C |x_1-x_2|^{\beta}$ for all $x_1,x_2 \in [0,1]$. 
In particular, we denote by $\mathcal{H}^{2+\beta,1+\beta/2}([0,1] \times [0,T))$, the space of maps $u : [0,1] \times [0,T) \to \mathbb{R}$ such that $u_{xx}$ and $u_t$ are Hölder continuous in $x$ with exponent $\beta$ uniformly in $t$ and Hölder continuous in $t$ with exponent $\beta/2$ uniformly in $x$.  
We denote by $\mathcal{H}_{loc}^{2+\beta,1+\beta/2}([0,1]\times [0,T))$ the space of maps $u :[0,1] \times [0,T) \to \mathbb{R}$ such that the restriction of $u$ to any bounded subset $[0,1] \times L \subset [0,1] \times [0,T)$ belongs to $\mathcal{H}^{2+\beta,1+\beta/2}([0,1] \times L)$. 

\section{Problem Formulation}\label{system description}
\subsection{Class of systems}\label{open_loop}
Consider the following quasilinear parabolic PDE 
\begin{empheq}[left=\Sigma:\ \left\{,right=\right.]{align}
&u_t   = \varepsilon(u)u_{xx} + \sum_{i=1}^{n}\gamma_i u^i u_x + u^p, \label{one}\\
&u_x(0)= v_0(u(0)), \ \ u_x(1) = v_1(u(1)), \label{two}\\
&u(x,0)= u_0(x), \label{three}
\end{empheq}
where $x\in [0,1]$ denotes the space variable, $t\geq 0$ the time variable, and $u\in \mathbb{R}$ the state variable. 

The term $\varepsilon(u)u_{xx}$ is referred to as diffusion, with $\varepsilon$ being the diffusion coefficient. The term $\sum_{i=0}^{n}\gamma_i u^i u_x$, for $n \in \mathbb{N}$ and $\gamma_i \in \mathbb{R}$ for all $i \in \{0,1,...,n\}$, constitutes the convection and $\{\gamma_i\}_{i=0}^{n}$ are the convection coefficients. This term models transport phenomena. The polynomial $u^p$, where $p \in (1,+\infty)$, is the reaction. Such a reaction term is employed, e.g., in modeling autocatalysis in chemical kinetics, where the presence of a substance amplifies its own production rate \cite{chemical1}. This term is responsible for finite-time blow-up phenomena, and is often referred to as superlinear, because 
\begin{align*}
\lim_{u\to +\infty} \frac{u^{p}}{u} = +\infty.
\end{align*}

Furthermore, \eqref{two} represents Neumann-type boundary conditions, with $v_1,v_o: \mathbb{R}\to \mathbb{R}$ some feedback laws to be designed.

Finally, \eqref{three} accounts for the initial condition $u_o\in  \mathcal{H}^{2+\beta}[0,1]$ for some $\beta \in (0,1)$. This regularity requirement is standard in the literature and is necessary to establish the existence of classical solutions \cite{lady}. 
\begin{definitionn}\label{def1}
A classical solution to $\Sigma$ defined on $[0,T)$, for some $T\in (0,+\infty]$, is any map $u \in \mathcal{H}^{2+\beta,1+\beta/2}_{loc}([0,1]\times [0,T))$ that verifies \eqref{one} for all $(x,t)\in (0,1)\times (0,T)$, \eqref{two} for all $t\in [0,T)$, and \eqref{three} for all $x\in [0,1]$. A classical solution is said to be complete if it is defined on $[0,+\infty)$.
\end{definitionn}

\begin{Remark}
According to the latter definition, a necessary condition for the existence of solutions is
\begin{align}
u_o'(0) = v_o(u_o(0)), \ \ u_o'(1) = v_1(u_o(1)). \label{compatibility}
\end{align}
We refer to \eqref{compatibility} as the compatibility condition of order $0$.
\end{Remark}

Our goal is to prevent, by an appropriate design of $v_o$ and $v_1$, and for a class of initial conditions to be specified later, finite-time blow up of the solutions to $\Sigma$. To state precisely this objective, we first need to clarify the meaning of finite-time blow up in this article, as different types of blow-up behavior can occur, as discussed in the remark below. A solution $u$ to $\Sigma$ is said to blow-up in finite time if there exists a time $0<T_m<+\infty$, called the blow-up time, such that 
\begin{align}
\lim_{t\to T_m^-}|u(\cdot,t)|_{\infty} = +\infty. \label{blow_up}
\end{align}
\begin{Remark}
If $\varepsilon$ belongs to $\mathcal{C}^2(\mathbb{R};\mathbb{R}_{>0})$ and is lower bounded by a positive constant, then the maximal time of existence of $u$, if it is finite, corresponds to the blow-up time \cite{noquenching1,noquenching2}. Namely, if $u$ is defined on a maximal interval $[0,T_m)$ for some $0<T_m<+\infty$, then \eqref{blow_up} holds. Not all quasilinear parabolic equations exhibit this property. Indeed, for some equations, one of the partial derivatives $u_x$, $u_{xx}$ or $u_t$ may blow-up in finite time while $|u|_{\infty}$ remains bounded. This phenomena is referred to as \textit{quenching} \cite{quenching}.
\end{Remark}

In general, there exists infinitely many sufficiently-'large' initial conditions $u_o$ such that, given any two maps $v_1$ and $v_o$, the corresponding solutions to $\Sigma$, if they exist, blow up in finite time. This holds true if, e.g., $\varepsilon$ is constant, $\gamma_i:=0$ for all $i\in \{0,1,...,n\}$, and $p:=2$; see \cite[Section 6.2]{volterra1}. For this reason, we do not seek global stability guarantees. Instead, as precisely stated in the next section, our objective is to design $v_1$ and $v_o$ to ensure stability properties for the origin $\{u=0\}$ for $\Sigma$ with an estimate of the region of attraction, while also guaranteeing $u(x,t)\geq 0$ for all $(x,t)\in [0,1]\times [0,+\infty)$ whenever $u_o(x)\geq 0$ for all $x\in [0,1]$. Such guarantees are still important since, when $(v_1,v_o):=0$, the solutions to $\Sigma$ can blow up in finite time regardless of the 'size' of $u_o$. This is best illustrated in the following example, where we consider
\begin{align}
u_t = \varepsilon(u) u_{xx} + u^2, \label{eg_illu}
\end{align}
subject to
\begin{align*}
u_x(0) = u_x(1) = 0.
\end{align*}
We assume that $\varepsilon$ is of class $\mathcal{C}^2$, lower-bounded by a positive constant, and $\varepsilon'(s)\leq 0$ for all $s\in \mathbb{R}$. Upon integration with respect to $x$, we obtain on the interval of existence of $u$,
\begin{align}
&\frac{d}{dt}\int_{0}^{1}u(x)dx
= \int_{0}^{1}\varepsilon(u(x))u_{xx}(x)dx + \int_{0}^{1}u(x)^2dx \nonumber \\
&= - \int_{0}^{1}\varepsilon'(u(x))u_x(x)^2dx+ \int_{0}^{1}u(x)^2dx.\nonumber 
\end{align}
As a result, using Cauchy-Schwarz inequality, we get 
\begin{align}
\frac{d}{dt}\int_{0}^{1}u(x)dx \geq \bigg( \int_{0}^{1}u(x)dx\bigg)^2. \label{blow_up_1}
\end{align}
In particular, if 
$$\int_{0}^{1}u_o(x)dx>0,$$ then $u$ blows up in finite time since $$\int_{0}^{1}u(x,t)dx\leq |u(\cdot,t)|_{\infty}.$$
In other words, for this particular example, it is sufficient to have the average of $u_o$ positive for finite-time blow up to occur. Furthermore, according to \eqref{blow_up_1}, the maximal time of existence of $u$, denoted $T_m$, verifies 
\begin{align}
T_m \leq \bigg(\int_{0}^{1}u_o(x)dx\bigg)^{-1}. \label{max_time}
\end{align}
Interestingly, $T_m$ is bounded by a constant that is independent of $\varepsilon$. Some of our results will be illustrated on \eqref{eg_illu} in Section \ref{simu}, for some state-dependent diffusion coefficient.

\subsection{Control-design problem}

The objective in this work is to design $v_1$ and $v_o$ to guarantee the following four properties.

\begin{Property}[Stability]\label{prop1}
There exist $\lambda_o,\lambda_1, \omega >0$ and $\mu \geq 0$, such that, if
\begin{align}
\kappa_o :=&~ \bigg[\int_{0}^{1}\big(u_o(x)^2+u_o'(x)^2\big)dx + \lambda_1u_o(1)^2+ \frac{\mu}{2}u_o(1)^4\nonumber \\
&~\quad + \lambda_o u_o(0)^2+ \frac{\mu}{2}u_o(0)^4\bigg]^{1/2} \leq \sqrt{2\omega}, \label{init_func}
\end{align}
and if $u$ is a complete solution to $\Sigma$, then there exists a constant $\sigma(\kappa_o)>0$ such that, 
\begin{align}
\int_{0}^{1}u(x,t)^2dx &\leq \bigg( \int_{0}^{1}u_o(x)^2dx\bigg)\exp^{-\sigma(\kappa_o) t} \quad \forall t\geq 0,  \label{decay_L2} \\
|u(\cdot,t)|_{\infty}^2 &\leq 2\kappa_o\left(\sqrt{\int_{0}^{1}u_o(x)^2dx}\right) \exp^{-(\sigma(\kappa_o)/2)t} \nonumber \\
&~+ \bigg(\int_{0}^{1}u_o(x)^2dx\bigg)\exp^{-\sigma(\kappa_o) t} \quad \forall t\geq 0.\label{decay_max}
\end{align}
\end{Property}

\begin{Property}[Well-posedness]\label{prop2}
Let $\kappa_o$ and $\omega$ be introduced in Property \ref{prop1}, such that $\kappa_o\leq \sqrt{2\omega}$ and the compatibility condition \eqref{compatibility} holds. Then, there exists a unique complete solution to $\Sigma$.
\end{Property}

\begin{Property}[Positivity]\label{prop3}
If $u_o(x)\geq 0$ for all $x\in [0,1]$ and $u$ is a complete solution to $\Sigma$, then $u(x,t)\geq 0$ for all $(x,t)\in [0,1]\times [0,+\infty)$.
\end{Property}

\begin{Property}[Invertibility]\label{prop4}
The maps $v_o$ and $v_1$ are invertible. That is, there exist maps $d_o:\mathbb{R}\to \mathbb{R}$ and $d_1:\mathbb{R}\to \mathbb{R}$ such that, given any solution $u$ to $\Sigma$,
\begin{equation*}
\begin{aligned}
&u_x(0) = v_o(u(0)) \Leftrightarrow u(0) = d_o(u_x(0)), \\
&\text{and} \ \ u_x(1) = v_1(u(1)) \Leftrightarrow u(1) = d_1(u_x(1)).
\end{aligned}
\end{equation*}
\end{Property}

Properties \ref{prop1} and \ref{prop2} indicate that if $u_o$ is sufficiently small, meaning that $\kappa_o\leq \sqrt{2\omega}$ - which constrains the $H^1$ norm of $u_o$ - then $\Sigma$ admits a unique complete solution under our boundary controllers. Along this solution, both the $L^2$ norm and the max norm remain bounded and decay exponentially to zero. It is important to note that boundedness of $t \mapsto \int_{0}^{1}u(x,t)^2dx$ does not necessarily imply the boundedness of $t\mapsto |u(\cdot,t)|_{\infty}$. However, Agmon’s inequality ensures that if both $t\mapsto \int_{0}^{1}u(x,t)^2dx$ and $t\mapsto \int_{0}^{1}u_x(x,t)^2dx$ are bounded, then $t\mapsto |u(\cdot,t)|_{\infty}$ is also bounded. This observation highlights the need to analyze $ \int_{0}^{1} u(x)^2 dx$ alongside $\int_{0}^{1} u_x(x)^2 dx$ when proving Property \ref{prop1}, providing insight into why $\kappa_o$ depends on $\int_{0}^{1} u_o'(x)^2 dx$.

Ensuring Property \ref{prop3} is essential, e.g., when negative values of the state are physically meaningless. This is the case, e.g., if $\Sigma$ is used to model autocatalysis in chemical kinetics \cite{chemical1,chemical2}.

Finally, Property \ref{prop4} means that the system $\Sigma'$, which is defined as $\Sigma$ with \eqref{two} replaced by one of the following boundary conditions 
\begin{align} 
&u(0) = d_o(u_x(0)), \quad u(1) = d_1(u_x(1)) \quad \forall t \geq 0, \label{dir1} \\
&u(0) = d_o(u_x(0)), \quad u_x(1) = v_1(u(1)) \quad \forall t \geq 0, \label{dir2}\\
&u_x(0) = v_o(u(0)), \quad u(1) = d_1(u_x(1)) \quad \forall t \geq 0, \label{dir3}
\end{align}
verifies Properties \ref{prop1}-\ref{prop3}. The motivation for ensuring the invertibility of $v_o$ and $v_1$ stems from practical considerations. For instance, in fluid flow control systems, Dirichlet-type actuation is often more practical to implement via blowing and suction of fluid than Neumann-type actuation \cite{dirichlet1,dirichlet2}, whereas in thermal systems, imposing a heat flux via a Neumann-type actuation is more practical than directly imposing the temperature (Dirichlet-type actuation) \cite{heat}; see also \cite[Equations (3.12) and (3.13)]{burgers}, \cite[Remark 2.3]{dec} and \cite[Equation (2.25)]{ks}.

To guarantee Properties \ref{prop1}-\ref{prop4}, we impose the following assumption on $\varepsilon$.

\begin{assumption}\label{ass1}
$\varepsilon$ is differentiable. Moreover, there exists a continuous non-decreasing map $\bar{\varepsilon}' : \mathbb{R}_{\geq 0}\to \mathbb{R}_{\geq 0}$ such that, for any $s\geq 0$,
\begin{align*}
\sup_{u\in [-s,s]}|\varepsilon'(u)|\leq \bar{\varepsilon}'(s).
\end{align*}
Finally, there exists a constant $\underline{\varepsilon}>0$ such that 
\begin{align*}
    \underline{\varepsilon} \leq \varepsilon (u) \quad \forall u\in \mathbb{R}.
\end{align*}
\end{assumption}
\begin{remark}
As we shall see later, the value of $\underline{\varepsilon}$ affects mainly the basin of attraction and the performance of our controllers. That is, the smaller $\underline{\varepsilon}$ is, the smaller the basin of attraction will be, and the slower the decay rate $\sigma(\kappa_o)$ in Property \ref{prop1} will be. Whereas, a larger $\underline{\varepsilon}$ leads to a larger basin of attraction and a faster decay rate.
\end{remark}

To establish Property \ref{prop2}, we consider, in addition to Assumption \ref{ass1}, the following regularity condition.
\begin{assumption}\label{ass2}
The second-order derivative $\varepsilon''$ exists and is continuous.
\end{assumption}

To ensure Property \ref{prop3}, we do not need Assumption \ref{ass2}. Instead, we require the following weaker condition. 
\begin{assumption}\label{ass1/2}
$\varepsilon$ is differentiable and $\varepsilon' \in \mathcal{H}_{loc}(\mathbb{R})$.
\end{assumption}

\section{Main results}\label{main results}
In this section, we present our main results. Due to space constraints, we provide only a sketch of the proof of Theorem \ref{thm1}. The full proof is available in \cite{journal_v}.

\subsection{Statement of the main results}
To state our results, we introduce the feedback laws $v_1$ and $v_o$, as well as $\kappa_o$ and $\omega$ in Property \ref{prop1}. That is, we let
\begin{align}
v_1(u(1)) :=&~ -\lambda_1u(1)-\mu u(1)^3, \label{v1} \\
v_o(u(0)) :=&~ \lambda_ou(0)+\mu u(0)^3, \label{v2}
\end{align}
where
\begin{align}
\lambda_1 &:= \frac{2(m_1+k_1)+|\gamma_1|}{2\underline{\varepsilon}}, \label{const1}\\
\lambda_o &:= \frac{2(m_o+k_o)+|\gamma_1|}{2\underline{\varepsilon}}, \label{const2}\\
\mu &:= \frac{1}{\underline{\varepsilon}}\bigg( \frac{|\gamma_1|}{2} +\sum_{i=2}^{n}\frac{|\gamma_i | M^{i-2}}{i+2}\bigg).\label{const3}
\end{align}
Here, $M, m_1, m_o>0$ and $k_1,k_o\geq 0$ are control parameters. Furthermore, $\omega >0 $ is any constant that verifies
\begin{align}
m_1+m_o & > 2(6\omega)^{(p-1)/2}+(6\omega)^{p-1}/\underline{\varepsilon}, \label{omega_1} \\
\underline{\varepsilon}-2(m_1+m_o)& > \sqrt{6\omega}\bar{\varepsilon}'(\sqrt{6\omega}) +\frac{n}{2\underline{\varepsilon}}\bigg(\sum_{i=1}^{n}\gamma_i^2(6\omega)^i\bigg), \label{omega2}\\
\sum_{i=2}^{n}\frac{|\gamma_i|}{i+2}M^{i-2} & > \sum_{i=2}^{n}\frac{|\gamma_i|}{i+2}(6\omega)^{(i-2)/2} \nonumber \\
&~\quad \text{if $\exists i\in \{2,3,...,n\}$ \text{s.t.}  $\gamma_i\neq 0$,} \label{omega3}
\end{align}
\begin{Remark}
Since $\bar{\varepsilon}'$ is non-decreasing, then the inequalities \eqref{omega_1}-\eqref{omega3} hold for $\omega$ sufficiently small, provided that the control gains $m_1$ and $m_o$ are chosen such that $\underline{\varepsilon}  > 2(m_1+m_2)$. E.g., one could select
\begin{align}
m_1 = m_o := \frac{\underline{\varepsilon}}{4+\delta} \quad \text{for some $\delta>0$}. \label{gains}
\end{align}
\end{Remark}
\begin{Theorem}\label{thm1}
Let Assumption \ref{ass1} be verified, and let $v_1$ and $v_o$ be given by \eqref{v1} and \eqref{v2}, respectively. Then, Property \ref{prop1} holds with $\lambda_1,\lambda_o,\mu$ given by \eqref{const1}, \eqref{const2}, and \eqref{const3} respectively, $\omega>0$ being any constant verifying \eqref{omega_1}-\eqref{omega3}, and 
\begin{align}
\sigma := m_1+m_o-2\left(\sqrt{3}\kappa_o\right)^{(p-1)} > 0, \label{sigma}
\end{align}
where $\kappa_o$ is defined in \eqref{init_func}.

Additionally, if Assumption \ref{ass2} holds, then Property \ref{prop2} (well-posedness) is also verified. Finally, if, instead of Assumption \ref{ass2}, we suppose only that Assumptions \ref{ass1} and \ref{ass1/2} are verified, then Property \ref{prop3} (positivity) holds. 
\end{Theorem}

\begin{Proposition}
Let Assumption \ref{ass1} holds, and $v_1$ and $v_o$ be defined in \eqref{v1} and \eqref{v2} respectively. Then, Property \ref{prop4} (Invertibility) is verified. In particular, 
\begin{itemize}
    \item If $\mu =0$, then 
    \begin{align}
d_1(u_x(1)) &:= -\frac{u_x(1)}{\lambda_1}, \label{lin1} \\
d_o(u_x(0)) &:= \frac{u_x(0)}{\lambda_o}. \label{lin2}
\end{align}
\item Otherwise, 
\begin{align}
d_1(u_x&(1)) := \sqrt[3]{-\frac{u_x(1)}{2\lambda_1}+\sqrt{\frac{u_x(1)^2}{4\lambda_1^2} + \frac{\lambda_1^3}{27\mu^3}}} \nonumber \\
&~+\sqrt[3]{-\frac{u_x(1)}{2\lambda_1}-\sqrt{\frac{u_x(1)^2}{4\lambda_1^2} + \frac{\lambda_1^3}{27\mu^3}}}. \label{cubicc1} 
\end{align} 
\begin{align}
d_o(u_x&(0)) := \sqrt[3]{\frac{u_x(0)}{2\lambda_o}+\sqrt{\frac{u_x(0)^2}{4\lambda_o^2} + \frac{\lambda_o^3}{27\mu^3}}} \nonumber \\
&~+\sqrt[3]{\frac{u_x(0)}{2\lambda_o}-\sqrt{\frac{u_x(0)^2}{4\lambda_o^2} + \frac{\lambda_o^3}{27\mu^3}}}. \label{cubicc2}
\end{align}
\end{itemize}
\end{Proposition}
\begin{proof}
If $\mu=0$, then
\begin{align}
v_1(u(1)) = -\lambda_1u(1) \quad \text{and} \quad v_o(u(0)) = \lambda_o u(0), \label{linear}
\end{align}
which implies \eqref{lin1}-\eqref{lin2}.

Now, suppose that $\mu \neq 0$. To find $d_1$ (resp. $d_2$), it is sufficient to rewrite $u(1)$ (resp. $u(0)$) as a function of $u_x(1)$ (resp. $u_x(0)$). To do so, note that we have 
\begin{align}
u(1)^3 + \frac{\lambda_1}{\mu}u(1) + \frac{1}{\lambda_1}u_x(1) = 0. \label{polynomial}
\end{align}
The latter is a cubic polynomial equation in $u(1)$. Its discriminant is given by
\begin{align}
\Delta := \frac{u_x(1)^2}{4\lambda_1^2} + \frac{\lambda_1^3}{27\mu^3} \geq 0. 
\end{align}
Since the discriminant is nonnegative, then \eqref{polynomial} admits a unique real solution, which is given by \eqref{cubicc1}. We prove similarly \eqref{cubicc2}.
\end{proof}

\subsection{Sketch of proof of Theorem \ref{thm1}}
\textit{Establishing Property \ref{prop1}.} The proof is divided into five steps. First, we consider the Lyapunov functional candidate
\begin{align}
V(u) := \frac{1}{2}\int_{0}^{1}u(x)^2dx. \label{V_def}
\end{align}
and we show that, along $u$, we have 
\begin{align*}
& \dot{V}  \leq -\bigg(m_1+m_o-2|u|_{\infty}^{p-1}\bigg)V -k_1 u(1)^2 - k_o u(0)^2
\\ &
- \bigg(\underline{\varepsilon}-|u|_{\infty}\bar{\varepsilon}(|u|_{\infty}) - 2(m_o+m_1)\bigg)\int_{0}^{1}u_x(x)^2dx \\
&+ u(1) \varepsilon(u(1)) \times\\
& 
\bigg[ v_1(u(1))  
+ \lambda_1 u(1) + \bigg( \frac{|\gamma_1|}{2 \underline{\varepsilon}} + \sum_{i=2}^{n}\frac{|\gamma_i | |u|_{\infty}^{i-2}}{\underline{\varepsilon}(i+2)}\bigg)u(1)^3\bigg]\\
& 
+ u(0)\varepsilon (u(0)) \times 
\\ \hspace{-0.3cm} & 
\bigg[ - v_o(u(0)) +\lambda_o u(0) + 
\bigg( \frac{|\gamma_1|}{2 \underline{\varepsilon}} + \sum_{i=2}^{n}\frac{|\gamma_i | |u|_{\infty}^{i-2}}{\underline{\varepsilon} (i+2)}\bigg)u(0)^3\bigg].
\end{align*}
As a result, under \eqref{v1} and \eqref{v2}, the following inequality holds
\begin{align}
\dot{V} \leq&~ -\bigg(m_1+m_o-2|u|_{\infty}^{p-1}\bigg)V - \bigg(\underline{\varepsilon}-|u|_{\infty}\bar{\varepsilon}'(|u|_{\infty}) \nonumber \\
&~-2(m_o+m_1)\bigg) \int_{0}^{1}u_x(x)^2dx \nonumber \\
&~-\varepsilon(u(1))\bigg[ \sum_{i=2}^{n}\frac{|\gamma_i|}{i+2}\bigg(M^{i-2}-|u|_{\infty}^{i-2}\bigg)\bigg] u(1)^4\nonumber \\
&~-\varepsilon(u(0))\bigg[ \sum_{i=2}^{n}\frac{|\gamma_i|}{i+2}\bigg(M^{i-2}-|u|_{\infty}^{i-2}\bigg)\bigg] u(0)^4 \nonumber \\
&~-k_1u(1)^2-k_ou(0)^2. \label{ineq_V_fin2}
\end{align}
In light of \eqref{ineq_V_fin2}, we need to analyze $V$ and $|u|_{\infty}$ simultaneously. To do so, we exploit Agmon's inequality, which states that $|u|_{\infty}$ can be upperbounded using $\int_{0}^{1}u(x)^2dx$ and $\int_{0}^{1}u_x(x)^2dx$. The second step of the proof consists then in deriving a differential inequality on  
\begin{align}
H(u) :=&~ \frac{1}{2}\int_{0}^{1}u_x(x)^2dx + \frac{\lambda_1}{2}u(1)^2+ \frac{\mu}{4}u(1)^4\nonumber \\
&~+ \frac{\lambda_o}{2}u(0)^2+ \frac{\mu}{4}u(0)^4. \nonumber
\end{align}
We prove that $H$ verifies
\begin{align}
\dot{H}\leq \frac{n}{2\underline{\varepsilon}}\bigg(\sum_{i=1}^{n}\gamma_i^2|u|_{\infty}^{2i}\bigg)\int_{0}^{1}u_x(x)^2dx + \frac{|u|_{\infty}^{2p-2}}{\underline{\varepsilon}}V. \label{ineq_H_to_use}
\end{align}
The third step of the proof combines \eqref{ineq_V_fin2} and \eqref{ineq_H_to_use}, and uses Agmon's inequality, to deduce boundedness of $|u|_{\infty}$ when $\kappa_o\leq \sqrt{2\omega}$. In particular, we let
\begin{align*}
E := V+H.
\end{align*}
According to Agmon's inequality, we have 
\begin{align}
|u|_{\infty}^2 \leq 2\sqrt{2V\int_{0}^{1}u_x(x)^2dx}+2V \leq 6E, \label{agmon_used}
\end{align}
which implies that if $E$ is bounded, then $|u|_{\infty}$ is also bounded. Now, using \eqref{ineq_V_fin2}, \eqref{ineq_H_to_use}, and \eqref{agmon_used}, we obtain
\begin{align}
&\dot{E} \leq - \alpha (E)V-\Lambda (E)\int_{0}^{1}u_x(x)^2dx\nonumber \\
&-\varepsilon(u(1))\Gamma(E)u(1)^4-\varepsilon (u(0))\Gamma(E)u(0)^4,\label{E_ineq_use}
\end{align}
where
\begin{align*}
\alpha(E) :=&~ m_1+m_o-2(6E)^{(p-1)/2} - (6E)^{p-1}/\underline{\varepsilon},\\
\Lambda (E) :=&~ \underline{\varepsilon}-\sqrt{6E}\bar{\varepsilon}'\big(\sqrt{6E}\big) -2(m_o+m_1) 
\\ &
- n \bigg(\sum_{i=1}^{n}\gamma_i^2(6E)^{i}\bigg) / (2\underline{\varepsilon}), \\
\Gamma(E) :=&~ \sum_{i=2}^{n}\frac{|\gamma_i|}{i+2}\bigg(M^{i-2}-(6E)^{(i-2)/2}\bigg).
\end{align*}
By exploiting \eqref{E_ineq_use} and the condition $\kappa_o\leq \sqrt{2\omega}$, while noting that $E(0)=(1/2)\kappa_o^2$, we can show that 
\begin{align}
E(t) \leq E(0) \quad \forall t\geq 0. \label{stab_E}
\end{align}
The fourth step of the proof exploits \eqref{ineq_V_fin2}, \eqref{agmon_used}, and \eqref{stab_E} to establish \eqref{decay_L2}. Indeed, we can show that 
\begin{align}
\dot{V} \leq - \sigma(\kappa_o)V-k_1u(1)^2 - k_o u(0)^2, \label{k1ko}
\end{align}
which allows us to conclude on \eqref{decay_L2}. Finally, the last step uses \eqref{decay_L2} and \eqref{agmon_used} to obtain \eqref{decay_max}.

\textit{Establishing Property \ref{prop2}.} One of the key ingredients is inequality \eqref{decay_max} involving the max norm of $u$. Indeed, by exploiting the latter, we can show that, when $\kappa_o\leq \sqrt{2\omega}$, the well-posedness of $\Sigma$ is equivalent to the well-posedness of a 'cutoff' version of $\Sigma$, denoted $\bar{\Sigma}$. Specifically, we let 
\begin{equation*}
\bar{\Sigma} \ : \ \left\lbrace
\begin{aligned}
&u_t = \bar{\varepsilon} (u)u_{xx} + \bar{\gamma}(u) u_x + \bar{f}(u) \ \ x\in (0,1), \ t>0, \\
&u_x(0) = \bar{v}_o(u(0)), \ \ u_x(1) = \bar{v}_1(u(1)) \quad t\geq 0, \\
&u(x,0) = u_o(x) \quad x\in [0,1],
\end{aligned}
\right.
\end{equation*}
where 
\begin{equation*}
\begin{aligned}
\bar{\varepsilon}(u) &:= \Phi(u)\varepsilon(u)+1-\Phi(u), \\
\bar{\gamma}(u) &:= \Phi(u)\bigg(\sum_{i=0}^{n}\gamma_i u^i\bigg), \ \ \bar{f}(u) := \Phi(u) u^p, \\
\bar{v}_1(u) &:= \Phi(u) v_1(u), \ \ \bar{v}_o(u) := \Phi(u)v_o(u).
\end{aligned}
\end{equation*}
Here, $\Phi:\mathbb{R}\to [0,1]$ if of class $\mathcal{C}^2$ and verifies\footnote{We can select, e.g., $\Phi(u) := \varphi_1(u-\mathcal{M})$, with 
$$ \varphi_1(u) := 1 -\frac{\varphi_2(u)}{\varphi_2(u) + \varphi_2(1-u)},$$
and
\begin{equation*}
\varphi_2(u) := 
\left\lbrace 
\begin{aligned}
&\exp^{-u^{-1}} \quad &\text{if $u>0$}, \\
&0 \quad &\text{if $u\leq 0$}.
\end{aligned}
\right.
\end{equation*}
$\qquad \qquad \qquad $ \begin{tikzpicture}[scale=0.7]
  \centering
  \begin{axis}[
    xlabel={$u$},
    ylabel={$\Phi(u)$},
    xmin=-0.2, xmax=2.2,
    ymin=-0.25, ymax=1.5,
    grid=both,
    width=7cm,
    height=5cm,
    samples=50,
    smooth,
    domain=0.001:0.999,
    axis lines=middle,
    xtick=\empty,
    ytick=\empty
    ]
    % Define M value as 1 for plotting purposes
    \def\M{0.5}

    % Plot the function Phi(u) = φ₁(u-M)
    \addplot[thick, black, domain={\M+0.001}:{\M+0.999}] {1-(exp(-1/(x-\M))/(exp(-1/(x-\M)) + exp(-1/(1-(x-\M)))))};
    
    % Extend the function with constant values
    \addplot[thick, black, domain=0:0.5] {1}; % Phi = 1 for u ≤ 0.5
    \addplot[thick, black, domain=1.5:2] {0}; % Phi = 0 for u ≥ 1.5

    % Add vertical dotted lines at x = 0.5 and x = 1.5
    \addplot[dotted, gray] coordinates {(0.5, -0.1) (0.5, 1.1)};
    \addplot[dotted, gray] coordinates {(1.5, -0.1) (1.5, 1.1)};

    % Add labels for M and M+1
    \node[anchor=north] at (axis cs: 0.5, -0.05) {$\mathcal{M}$};
    \node[anchor=north] at (axis cs: 1.5, -0.05) {$\mathcal{M}+1$};
    \node[anchor=north] at (axis cs: -0.1, 1.1) {$1$};
  \end{axis}
\end{tikzpicture}
} 
\begin{equation} \label{cutoff}
\Phi (u) := 
\left\lbrace 
\begin{aligned}
&1 \quad \text{if $|u|\leq \mathcal{M}$}, \\
&0 \quad \text{if $|u|\geq \mathcal{M}+1$},
\end{aligned}
\right. 
\end{equation}
where 
\begin{align}
\mathcal{M} := \sqrt{2\sqrt{V(u_o)E(u_o)}+2V(u_o)} + \delta_{\mathcal{M}}, \label{cutoff_constant}
\end{align}
for some constant $\delta_{\mathcal{M}}>0$. We first establish the well-posedness of $\bar{\Sigma}$. Then, exploiting \eqref{decay_max}, we conclude on the well-posedness of $\Sigma$. A crucial step is to prove, using \eqref{decay_max}, that the unique complete solution $u$ to $\bar{\Sigma}$ verifies 
\begin{align}
|u(x,t)| \leq \mathcal{M} \quad \forall x\in [0,1], \ \forall t\geq 0. \label{bar_M}
\end{align}
Indeed, assuming we have established \eqref{bar_M}, it follows that
$$ \Phi(u(x,t))=1 \quad \forall x\in [0,1], \ \forall t\geq 0.$$
Consequently, we have 
\begin{equation}\label{all_terms}
\begin{aligned}
&\bar{\varepsilon}(u)=\varepsilon(u), \ \ \bar{\gamma}(u) = \sum_{i=0}^{n}\gamma_i u^i, \ \ \bar{f}(u) = u^p, \\
&\bar{v}_1(u(1))=v_1(u(1)), \ \ \bar{v}_o(u(0)) = v_o(u(0)),
\end{aligned}
\end{equation}
which implies that $u$ is also a complete solution to $\Sigma$.

For uniqueness, we note that according to Theorem \ref{thm1}, if $u$ is a complete solution to $\Sigma$, then \eqref{bar_M} holds since \eqref{decay_max} is verified on $[0,+\infty)$. Therefore, the identities in \eqref{all_terms} are satisfied, and $u$ must be a complete solution to $\bar{\Sigma}$. Since $\bar{\Sigma}$ admits a unique complete solution, we conclude that there exists at most one complete solution to $\Sigma$.

\textit{Establishing Property \ref{prop3}.}
The key is to observe that the maps $u\mapsto u^p$, $u\mapsto -v_1(u)$ and $u\mapsto v_o(u)$ are non-decreasing on $[0,+\infty)$. This property enables us to apply the parabolic comparison theorem \cite[Condition (H$_1$)]{maximum_principle} to complete the proof.

\section{Discussion}

\textit{Region of attraction and decay rate.} Inequalities \eqref{omega_1}-\eqref{omega3} show how the range of values of $\omega$, that ensure Property \ref{prop1}, depends on the control gains and the parameters of $\Sigma$. In particular, suppose, e.g., that $m_1$ and $m_o$ are given by \eqref{gains} and that $M:=M(\underline{\varepsilon})$ with $\lim_{\underline{\varepsilon}\to +\infty}M(\underline{\varepsilon})=+\infty$. Moreover, suppose that $\bar{\varepsilon}'(s)\leq \delta'$ for some $\delta'>0$ and for all $s \geq 0$. Given $p\in (1,+\infty)$, if we denote by $\Omega (\underline{\varepsilon})$ the supremum of the values of $\omega$ that ensure \eqref{omega_1}-\eqref{omega3}, then
$$\lim_{\underline{\varepsilon}\to +\infty} \Omega(\underline{\varepsilon})=+\infty.$$
This result contrasts with the case $(v_1,v_o):=0$, where finite-time blow up occurs independently of how large $\underline{\varepsilon}$ and how small $\kappa_o$ may be; see Section \ref{open_loop}. Furthermore, the expression for the decay rate $\sigma(\kappa_o)$ given in \eqref{sigma} implies that
$$ \sigma(\kappa_o) \to +\infty \quad \text{as} \quad \underline{\varepsilon}\to +\infty.$$

\textit{Linear vs nonlinear feedback.} The feedback laws $v_1$ and $v_o$ are nonlinear (i.e., $\mu \neq 0$) only if $\Sigma$ contains nonlinear convective terms (i.e., $\gamma_i\neq 0$ for some $i\in \{1,2,...,n\}$). If $\gamma_i=0$ for all $i\in \{1,2,...,n\}$, then the equations $u_x(1)=v_1(u(1))$ and $u_x(0)=v_o(u(0))$ reduce to the so-called Robin boundary conditions
\begin{align*}
u_x(1) = - \lambda_1 u(1), \quad u_x(0)=\lambda_o u(0).
\end{align*}

\textit{Stability enhancement via $k_o$ and $k_1$}. The importance of the gains $k_o$ and $k_1$ can be deduced from \eqref{k1ko}. According to the latter, to ensure \eqref{decay_L2}, it is sufficient to set $k_1=k_o=0$. However, we claim that in some cases, the inequality \eqref{decay_L2} can be strict if we let $k_1>0$ or $k_o>0$. Indeed, by integrating \eqref{k1ko}, we get, for any $t> 0$, 
\begin{align*}
V(t) \leq&~ V(0)\exp^{-\sigma(\kappa_o)t}-k_1\int_{0}^{t}\exp^{-\sigma(\kappa_o)(t-\tau)}u(1,\tau)^2d\tau \\
&~-k_o\int_{0}^{t}\exp^{-\sigma(\kappa_o)(t-\tau)}u(0,\tau)^2d\tau.
\end{align*}
Suppose that there exists a subset $I_t$ of $(0,t)$ of non-zero measure, such that, for each $\tau\in I_t$, $u(1,\tau)\neq 0$ or $u(0,\tau)\neq 0$. It implies that 
\begin{align*}
0 >&~-k_1\int_{0}^{t}\exp^{-\sigma(\kappa_o)(t-\tau)}u(1,\tau)^2d\tau \\
&~-k_o\int_{0}^{t}\exp^{-\sigma(\kappa_o)(t-\tau)}u(0,\tau)^2d\tau \nonumber,
\end{align*}
which further implies that 
$$V(t) < V(0)\exp^{-\sigma(\kappa_o)t}.$$
The existence of a time $t>0$ and a subset $I_t$ of $(0,t)$ as above is guaranteed under \eqref{v1}-\eqref{v2}, for an initial condition $u_o$ such that $\kappa_o\leq \sqrt{2\omega}$, if the corresponding solution to $\Sigma$ blows-up in finite time when the boundary inputs are set to zero for almost all time. This is holds true if we take, e.g., $p=2$, $\gamma_i:=0$ for all $i\in \{1,2,...,n\}$, $\varepsilon'\leq 0$, and $u_o$ such that $\int_{0}^{1}u_o(x)dx>0$ and $\kappa_o\leq \sqrt{2\omega}$ (see Section \ref{open_loop}, and note that finite-time blow up occurs if \eqref{blow_up_1} holds for almost all time on the interval of existence of $u$).

In the simulations results in Section \ref{simu}, we will see that, indeed, by increasing $k_1$ and $k_o$, the decay of $V$ and $|u|_{\infty}$ will become faster.

\textit{On the use of a single feedback controller.} Property \ref{prop1}, in many cases, can be guaranteed using a single feedback controller. For instance, if $\gamma_{2k+1}:=0$ and $\gamma_{2k}\geq 0$ (resp. $\gamma_{2k}\leq 0$) for all $k\in \mathbb{N}$ such that $2k+1\in \{0,1,...,n\}$, then Property \ref{prop1} can be ensured by designing $v_1$ as in \eqref{v1} (resp., $v_o$ as in \eqref{v2}) while setting $v_o:=0$ (resp., $v_1:=0$). Moreover, in these cases, Properties \ref{prop2} and \ref{prop3} can also be guaranteed.

\section{Simulation Results}\label{simu}
In this section, we illustrate our results via numerical simulations on MATLAB. 

We discretize uniformly the interval $[0,1]$ into a finite set of points, with a spatial step $\Delta x := 0.01$. At the interior points, we use a central difference method to approximate $u_{xx}$, and an Euler backward method to approximate $u_x$. Furthermore, $u_x(0)$ (resp. $u_x(1)$) is approximated using an Euler forward (resp. backward) method. Finally, the resulting semi-discretized system is integrated using the ode15s solver. 

We consider the quasilinear parabolic PDE
\begin{align*}
u_t = \bigg(\delta_1+\delta_2 \exp^{-u^3}\bigg)u_{xx}+\gamma_1 uu_x + u^2,
\end{align*}
where $\delta_1>0$ and $\delta_2\geq 0$, and $\gamma_1\in \mathbb{R}$. The diffusion coefficient verifies Assumptions \ref{ass1}-\ref{ass1/2} with $\bar{\varepsilon}'(s):=2\delta_2s^2\exp^{s^3}$ and $\underline{\varepsilon}=\delta_1$. 

Furthermore, we let the initial condition 
\begin{align*}
u_o(x) := \delta_o \big( 1-\cos(2\pi x)\big) \quad \forall x\in [0,1],
\end{align*}
for some $\delta_o >0$.

\begin{figure}
    \centering
    \includegraphics[width=1\linewidth]{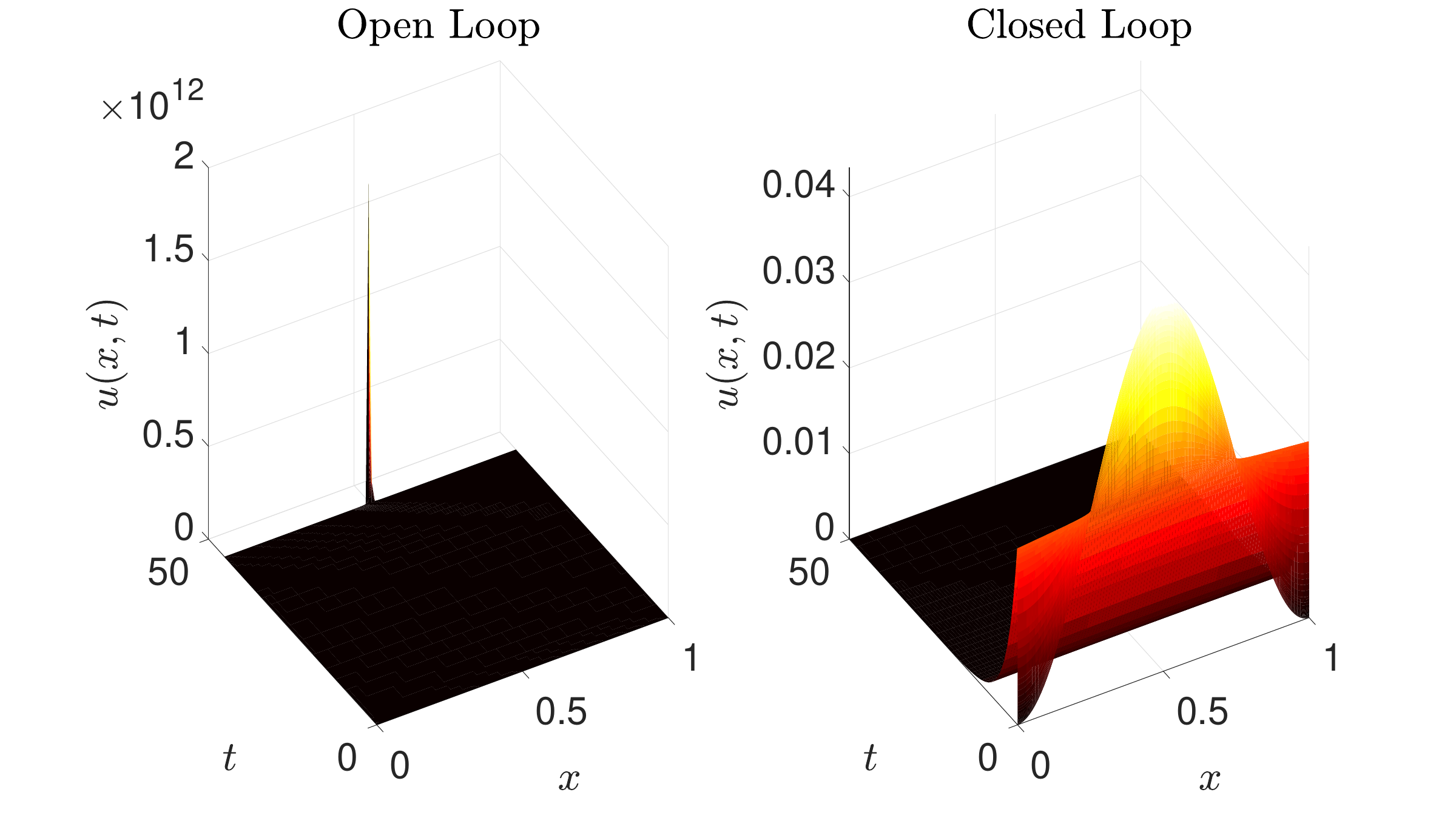}
    \caption{Left: response of $\Sigma$ to $u_x(1)=u_x(0):=0$, with $\gamma_1:=0$. Right: response of $\Sigma$ to $u_x(1):=-\lambda_1 u(1)$ and $u_x(0):=\lambda_o u(0)$, with $\gamma_1:=0$.}
    \label{fig1}
\end{figure}

First, we suppose that $\gamma_1:=0$, and we set $(v_1,v_o):=0$. In this case, since 
$$\int_{0}^{1}u_o(x)dx = \delta_o >0, $$ 
then, as discussed in Section \ref{open_loop}, the state $u$ must blow up in finite time. This is illustrated in Figure \ref{fig1} (left), where we have chosen
\begin{align}
(\delta_o,\delta_1,\delta_2) := \left(\sqrt{\frac{5 \times 10^{-3}}{\pi^2+3/4}},2,\frac{1}{24}\right). \label{deltas}
\end{align}
Moreover, note that in simulation, the maximal time of existence of $u$, denoted by $T_m$, is bounded by $50$. This is in agreement with \eqref{max_time}, which can be rewritten as
\begin{align*}
T_m \leq \bigg( \int_{0}^{1}u_o(x)dx\bigg)^{-1} &= \frac{1}{\delta_o} \approx 46.086.
\end{align*}

To prevent finite-time blow up, we employ the feedback controllers $v_1$ and $v_o$ in \eqref{v1} and \eqref{v2}, respectively. Since $\gamma_1=0$, then $\mu = 0$ and thus $v_1$ and $v_o$ reduce to the linear controllers 
\begin{align*}
&v_1(u(1)) := -\bigg(\frac{m_1+k_1}{\delta_1}\bigg) u(1), \\
&\text{and} \ \ v_o(u(0)) := \bigg(\frac{m_o+k_o}{\delta_1}\bigg) u(0), 
\end{align*}
where $m_1,m_o>0$ and $k_1,k_o\geq 0$ are control gains. Note that the compatibility condition \eqref{compatibility} holds, since $v_1(0)=v_o(0)=0$, and
\begin{align}
u_o(0)=u_o(1)=u_{o}'(0)=u_{o}'(1)=0. \label{compat_simu}
\end{align}
We choose $m_1$ and $m_o$ as in \eqref{gains}, with $\delta = 2$, i.e.,
\begin{align}
m_1 = m_o := \frac{\delta_1}{6}. \label{m1m2}
\end{align}
Furthermore, $\kappa_o$ in Theorem \ref{thm1} is given by 
\begin{align}
\kappa_o &:= \sqrt{\int_{0}^{1}\big(u_o(x)^2+u_o'(x)^2\big)dx} \nonumber \\
&= \delta_o \sqrt{\int_{0}^{1}\bigg[\big( 1-\cos(2\pi x)\big)^2+4\pi^2\sin(2\pi x)^2\bigg]dx} \nonumber \\
&= \delta_o\sqrt{2\pi^2+\frac{3}{2}}. \label{kappa_o_simu} 
\end{align}
According to Theorem \ref{thm1}, when $\kappa_o\leq \sqrt{2\omega}$, there exists a unique complete solution $u$ to $\Sigma$, along which, the inequalities \eqref{decay_L2} and \eqref{decay_max} hold, with the decay rate
\begin{align*}
\sigma := \frac{\delta_1}{3}-2\sqrt{3}\left(\sqrt{2\pi^2+\frac{3}{2}}\right)\delta_o.
\end{align*}
Here, $\omega>0$ is a constant verifying
\begin{align}
\frac{\delta_1}{3} &> 2\sqrt{6\omega} + \frac{6}{\delta_1}\omega, \label{sigma_pos}\\
\frac{\delta_1}{3} &>2\delta_2(6\omega)^{\frac{3}{2}}\exp^{(6\omega)^{3/2}}. \label{sigma_pos1/2} 
\end{align}
The latter two inequalities are verified for sufficiently small values of $\omega$, and \eqref{sigma_pos} ensures that $\sigma(\kappa_o)>0$. Moreover, note that supremum of the values of $\omega$ for which the latter two inequalities hold, goes to $+\infty$ as $\delta_1\to +\infty$.

We let $(\delta_o,\delta_1,\delta_2)$ be given by \eqref{deltas}, and note that \eqref{sigma_pos} and \eqref{sigma_pos1/2} are verified with 
$(1/2)\kappa_o^2 = \omega = 5\times 10^{-3}$. Furthermore, we let $(k_1,k_o):=0$, and plot the state $u$ in Figure \ref{fig1} (right). As expected, finite-time blow up is prevented and the state goes to zero. Moreover, we have $u(x,t)\geq 0$ for all $x\in [0,1]$ and all $t\geq 0$, as stated in Theorem \ref{thm1}.

\begin{figure}
    \centering
    \includegraphics[width=1\linewidth]{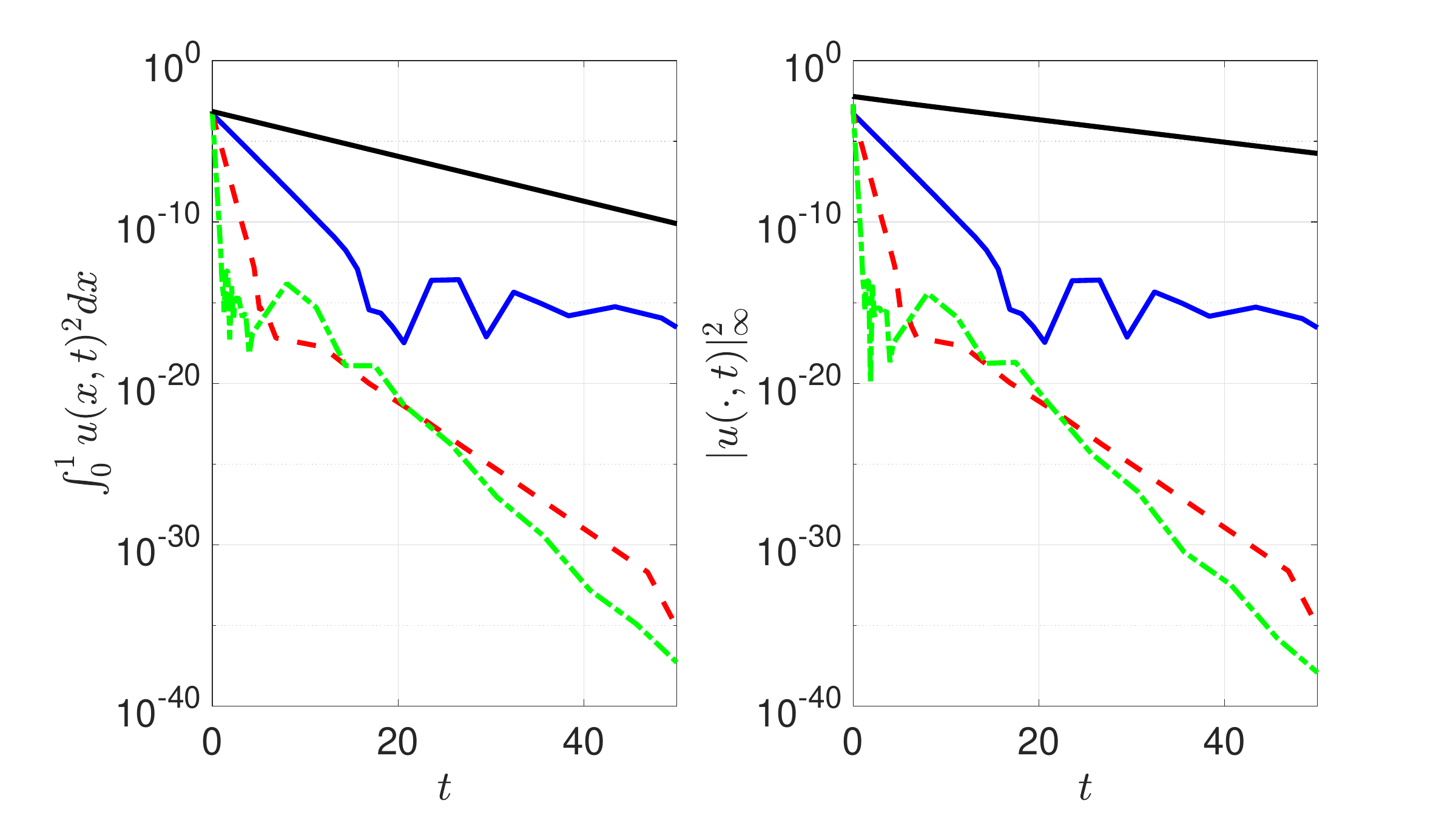}
    \caption{Left: theoretical bound on $t\mapsto \int_{0}^{1}u(x)^2dx$, i.e. the right-hand side of \eqref{decay_L2} (black); the map $t\mapsto \int_{0}^{1}u(x,t)^2dx$ for $k_1=k_o:=0$ (blue), $k_1=k_o:=1$ (red), and $k_1=k_o:=10$ (green). Right: theoretical bound on $t\mapsto |u(\cdot,t)|_{\infty}^2$, i.e. the right-hand side of \eqref{decay_max} (black); the map $t\mapsto |u(\cdot,t)|_{\infty}^2$ for $k_1=k_o:=0$ (blue), $k_1=k_o:=1$ (red), and $k_1=k_o:=10$ (green).}
    \label{fig2}
\end{figure}

Additionally, we plot in the left (resp., right) of Figure \ref{fig2}, in a semi-log scale, the map $t\mapsto \int_{0}^{1}u(x,t)^2dx$ (resp., $t\mapsto |u(\cdot,t)|_{\infty}^2$) as well as the right-hand side of inequality \eqref{decay_L2} (resp., \eqref{decay_max}) under the proposed boundary controllers, and for different values of the control gains $k_1$ and $k_o$, namely, for $k_1=k_o:=0$, $k_1=k_o:=1$, and $k_1=k_o:=10$. The simulations confirm that the inequalities \eqref{decay_L2} and \eqref{decay_max} are verified. Moreover, we note that larger values for $k_1$ and $k_o$ lead to a faster decay of $t\mapsto \int_{0}^{1}u(x,t)^2d$ and $t\mapsto |u(\cdot,t)|_{\infty}^2$.

\begin{figure}
    \centering
    \includegraphics[width=1\linewidth]{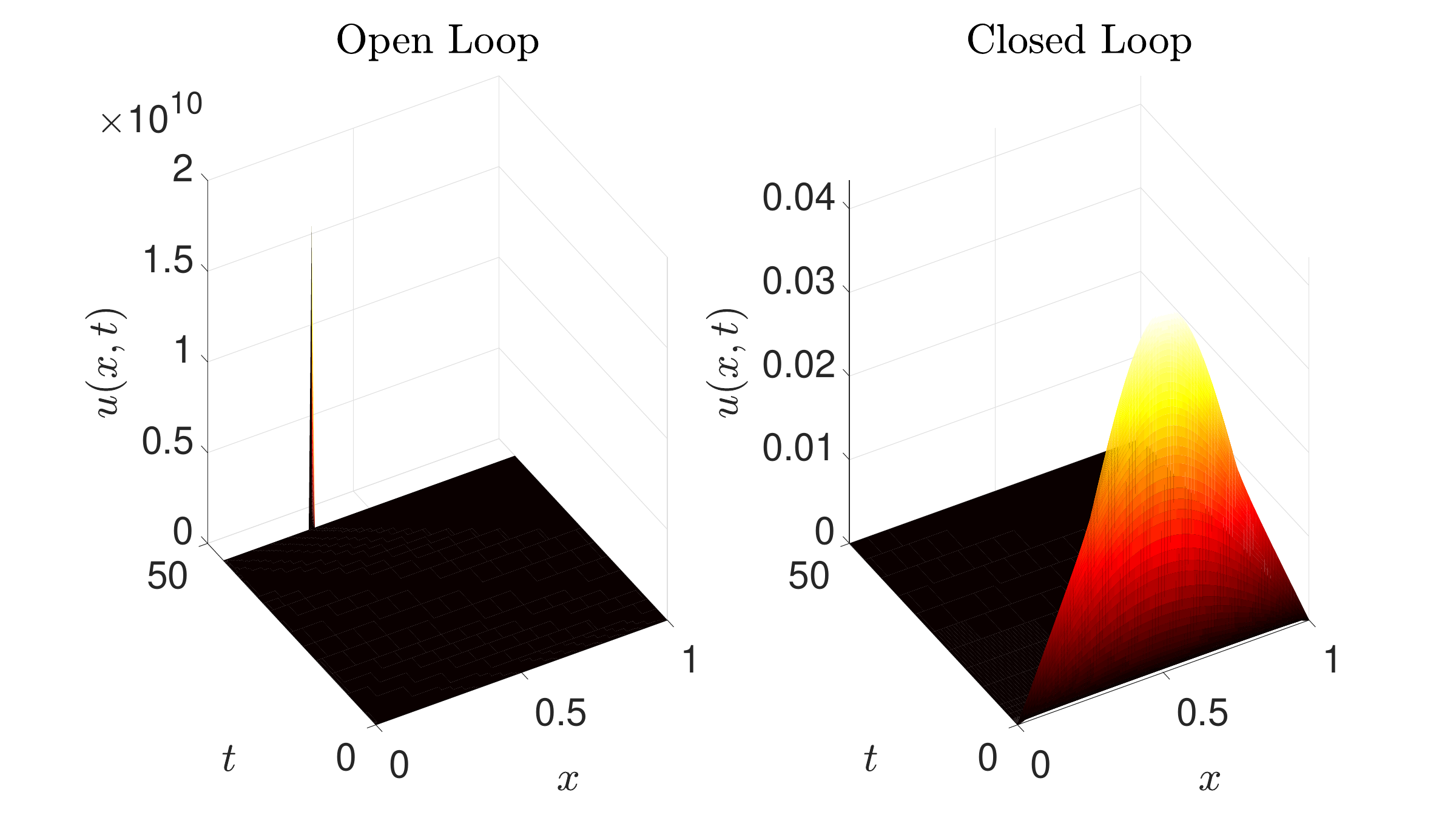}
    \caption{Left: response of $\Sigma$ to $u_x(1)=u_x(0):=0$, with $\gamma_1:=1$. Right: response of $\Sigma$ to $u_x(1):=-\lambda_1 u(1)-\mu u(1)^3$ and $u_x(0):=\lambda_o u(0)+\mu u(0)^3$, with $\gamma_1:=1$.}
    \label{fig3}
\end{figure}

We suppose now that $\gamma_1:=1$, i.e. the nonlinear convective term $uu_x$ is considered. In Figure \ref{fig3} (left), we simulate the response of $\Sigma$ under $(v_1,v_o):=0$, with $(\delta_o,\delta_1,\delta_2)$ as in \eqref{deltas}. We can see the finite-time blow up of $u$. To prevent it, we use the feedback controllers $v_1$ and $v_o$ in \eqref{v1} and \eqref{v2}, respectively. Since $\gamma_1=1$, then in this case, $v_1$ and $v_o$ are not linear, but instead given by 
\begin{align*}
v_1(u(1)) &:= -\bigg(\frac{2(m_1+k_1)+1}{2\delta_1}\bigg)u(1)-\frac{1}{2\delta_1}u(1)^3, \\
v_o(u(0)) &:= \bigg(\frac{2(m_o+k_o)+1}{2\delta_1}\bigg)u(0)+\frac{1}{2\delta_1}u(0)^3.
\end{align*}
Note that, since $v_1(0)=v_o(0)=0$, then the compatibility condition \eqref{compatibility} holds because of \eqref{compat_simu}.

We select $m_1$ and $m_o$ according to \eqref{m1m2}. Moreover, $\kappa_o$ is still given by \eqref{kappa_o_simu}. However, in this case, since $\gamma_1=1$, then $\omega$ needs to verify, in addition to \eqref{sigma_pos}, the following condition
\begin{align}
\frac{\delta_1}{3} >2\delta_2(6\omega)^{\frac{3}{2}}\exp^{(6\omega)^{3/2}}+ \frac{3\omega }{\delta_1}. \label{sigma_pos3}
\end{align}
The condition $\kappa_o\leq \sqrt{2\omega}$, for $\omega$ verifying \eqref{sigma_pos} and \eqref{sigma_pos3}, is satisfied for $(1/2)\kappa_o^2=\omega = 5\times 10^{-3}$ and $(\delta_o,\delta_1,\delta_2)$ as in \eqref{deltas}. The response of $\Sigma$ under our boundary controllers, with $k_1=k_o:=0$, is shown in Figure \ref{fig3} (right). This simulation result confirms Theorems \ref{thm1}, since $u$ converges to zero and $u(x,t)\geq 0$ for all $x\in [0,1]$ and all $t\geq 0$.

\section{Conclusion and Research Perspectives}
We studied stabilization of quasilinear parabolic PDEs that exhibit finite-time blow up in open loop. We designed Neumann-type boundary controllers as cubic polynomials of boundary measurements, ensuring $L^2$ exponential stability of the origin with an estimate of the region of attraction, boundedness and exponential decay of the state's max norm to zero, well-posedness, and positivity of solutions. Future research includes extending our approach to cases where $\varepsilon$ depends on $u_x$, generalizing to higher-dimensional domains, designing adaptive controllers when the coefficients are unknown, and considering in-domain measurement/actuation.

\end{document}